\documentclass{amsart}

\usepackage{amssymb}

\theoremstyle{plain}

\numberwithin{equation}{section}

\newtheorem*{teo}{Theorem A}
\newtheorem*{teob}{Theorem B}

\newtheorem{lem}[equation]{Lemma}

\newcommand{\Irr}{\operatorname{Irr}}
\newcommand{\Lin}{\operatorname{Lin}}

\newcommand{\Ker}{\operatorname{Ker}}

\theoremstyle{definition}

\newtheorem{example}[equation]{Example}
\begin{document}	
\title{Restriction of characters and products of characters}

\author{Edith Adan-Bante}

\address{Department of Mathematical Science, Northern Illinois University, 
 Watson Hall 320
DeKalb, IL 60115-2888, USA} 

\begin{abstract} Let $G$ be a finite $p$-group, for some prime $p$, 
 and $\psi, \theta \in \Irr(G)$ be irreducible
complex characters of $G$. It has been proved that if, in addition, $\psi,\theta$ are
faithful characters, then the product $\psi\theta$ is a multiple of an irreducible
or it is the nontrivial linear combination of at least $\frac{p+1}{2}$ distinct irreducible
characters of $G$. We show that if we do not require the characters to be faithful,
then given any integer $k>0$, we can always find a $p$-group $G$ and irreducible
characters $\Psi$ and $\Theta$ such that $\Psi\Theta$ is the nontrivial 
combination of exactly $k$ distinct irreducible characters. We do this
by translating examples of decompositions of restrictions of
characters into decompositions of products of characters.
\end{abstract}
\email{EdithAdan@illinoisalumni.org}

\subjclass{20c15}

\date{2007}
\maketitle

\begin{section}{Introduction}
Let $G$ be a finite group. Denote by $\Irr(G)$ the set of irreducible
complex 
characters of $G$. Let $\psi$ and $\theta$ be characters in $\Irr(G)$. 
Since a product of characters is a character, the product
$\psi\theta$, where $\psi\theta(g)=\psi(g)\theta(g)$
for all $g\in G$,
 is a character.
 Then the decomposition of
the character $\psi\theta$ into its distinct 
irreducible 
constituents $ \phi_1 ,\ \phi_2, \ldots, \phi_k, \in \Irr(G)$
has the form
\begin{equation*}	
 \psi\theta=  \sum_{i=1}^k n_i \phi_i
\end{equation*}
\noindent where $k>0$ and 
$n_i=(\psi\theta, \phi_i)>0$ is the multiplicity of 
$\phi_i$ in $\psi\theta$ for each $i=1, \ldots, k$. 
Let  $\eta(\psi\theta)= k $ be the
number of distinct 
 irreducible constituents 
of  the character $\psi\theta$. 

Given any subgroup $H$ of $G$, we denote by 
$\chi_H$ the restriction of the character $\chi$ of $G$ to $H$. 

Let $G$ be a finite $p$-group, where $p>2$ is a prime number. 
Given any two faithful characters $\psi,\theta \in \Irr(G)$, in Theorem A of 
\cite{edith1}
is proved that the product $\psi\theta$ is either a multiple of an irreducible, i.e
$\eta(\psi\theta)=1$, or $\psi\theta$ is the linear combination of at least
$\frac{p+1}{2}$ distinct irreducible constituents, i.e. 
$\eta(\psi\theta)\geq \frac{p+1}{2}$. Thus there is not $5$-group $G$ with 
faithful characters $\psi, \theta\in \Irr(G)$ such that 
$\psi\theta$ has exactly two distinct irreducible constituents. But can we
find a 5-group $P$ with characters
$\psi,\theta\in \Irr(P)$ such that 
$\psi\theta$ has exactly two distinct irreducible constituents? The answer is yes.
Moreover 

\begin{teo}
Fix  a prime number $p$. Let $k>0$ and $n_i>0$, $i=1,\ldots, k-1$, be integers.
Choose $r>0$ such that $p^r> \sum_{i=1}^{k-1} n_i$. Let
$n_k=p^r-\sum_{i=1}^{k-1} n_i$. Then there exists a finite $p$-group 
$G$ with characters $\Psi, \Theta \in \Irr(G)$ such that

$$\Psi\Theta = \sum_{i = 1}^k n_i\Phi_i, $$
 
\noindent where $\Phi_1, \Phi_2, \dots, \Phi_k$ are distinct irreducible
characters in $\Irr(G)$ and $n_i=(\Psi\Theta, \Phi_i)$.  Thus $\eta(\Psi\Theta)=k$.
\end{teo}

The main key for the previous result  is the following 

\begin{teob}
Let $P$ be a finite $p$-group, $Q<P$ be a subgroup of $P$ and $\psi\in \Irr(G)$. 
Assume that 

(i) $\psi_Q = \sum_{i = 1}^k n_i\phi_i, $

\noindent where $\phi_1, \phi_2, \dots, \phi_k$ are distinct irreducible
characters in $\Irr(Q)$ and $n_i = (\psi_Q, \phi_i) > 0$ for each $i =
1,2,\dots, k$. Then there exists a $p$-group $G$ with characters
$\Psi,\Theta \in \Irr(G)$ such that 

 (ii) $\Psi\Theta = \sum_{i = 1}^k n_i\Phi_i, $
 
\noindent where $\Phi_1, \Phi_2, \dots, \Phi_k$ are distinct irreducible
characters in $\Irr(G)$ and $n_i = (\Psi\Theta, \Phi_i)>0$ for each $i =
1,2,\dots, k$. 

\end{teob}
The previous result allows us to
translate all examples of decompositions of restrictions of
characters into decompositions of products of characters. Since
 the product $\psi\theta$ of two characters $\psi,\theta$ of a
group $G$ can be regarded as the restriction of the character $\psi
\times \theta$ of the direct product group $G \times G$ to the diagonal
subgroup $D(G) = \{ (g,g) | g \in G \}$, which is
another copy of $G$ inside $G \times G$, the converse also holds.

{\bf Acknowledgment.} 
Professor Everett C. Dade, in a personal communication, provided Theorem B. 
I thank him for that and for his 
suggestions.
\end{section}

\begin{section}{Proofs}
\begin{lem}
Given any finite $p$-group $P$ and any subgroup $Q$ of $P$, there
exists an elementary abelian $p$-group $E$ and $\lambda\in \Lin(E)$ such
that $P$  acts on $E$ as automorphisms and $Q$ is the stabilizer of $\lambda$ 
in $P$.  
\end{lem}
\begin{proof}
Let $S$ be a finite set on which $P$ acts transitively
with $Q$ as the stabilizer of one point $s \in S$. Then extend the
action of $P$ on $S$ to one of $P$ as automorphisms of the elementary
abelian $p$-group $E = <S>$ with $S$ as a basis. 
Finally, take
$\lambda$ to be any linear character on $E$ sending $s$ to some
primitive $p$-th root of unity, and every $s' \in S - \{s\}$ to $1$.
\end{proof}

{\bf Notation.} We will denote by $1_G$ the trivial character of $G$. We will also
denote by $\Lin(G)$ the set of linear characters of $G$. Also 
we denote by $1\in G$ the identity element of $G$.

\begin{proof}[Proof of Theorem B]
Let  $E$ be an elementary 
abelian $p$-group and $\lambda\in \Lin(E)$ such that
$P$ acts on $E$ as automorphisms of $E$ and the stabilizer $P_{\lambda}$ of
$\lambda$ in $P$ is $Q$. 
 We define $G$ to be the semi-direct product $P
\ltimes E$. We identify both $P$ and $E$ with their images in $G$. So
$$ G = PE = P \ltimes E.$$

We define $\Psi$ to be the inflation of $\psi \in \Irr(P)$ to an
irreducible character of $G$ with the values
\begin{equation}\label{Psi}
 \Psi(\sigma\tau) = \psi(\sigma) 
\end{equation}
\noindent  for any $\sigma \in P$ and $\tau \in E$. 

Because $Q$ fixes $\lambda$, the factor group $QE/\Ker(\lambda)$
is naturally isomorphic to $Q \times (E/\Ker(\lambda))$. It follows
that there is a unique character $\phi_i\lambda \in \Irr(QE)$, for each
$i = 1,2,\dots,k$, such that

 $$(\phi_i\lambda)(\sigma\tau) = \phi_i(\sigma)\lambda(\tau) $$

\noindent for any $\sigma \in Q$ and $\tau \in E$. The characters
$\phi_i\lambda$, for $i = 1,2,\dots,k$, are clearly distinct and lie
over $\lambda$. Because $QE = P_{\lambda}E$ is the stabilizer
$G_{\lambda}$ of $\lambda$ in $G$, Clifford theory tells us that the
characters

 $$ \Phi_i = (\phi_i\lambda)^G $$

\noindent induced by these $\phi_i\lambda$ are distinct and lie in $\Irr(G)$.

The trivial character $1_Q$ on $Q$ also defines an
irreducible character $1_Q\lambda$ of $QE$, with the values

 $$(1_Q\lambda)(\sigma\tau) = \lambda(\tau) $$

\noindent for any $\sigma \in Q$ and $\tau \in E$. This also induces an
irreducible character
\begin{equation}\label{Theta}
 \Theta = ( 1_Q\lambda)^G \in \Irr(G).
 \end{equation}

We now prove  that

 $$\Psi\Theta = \sum_{i = 1}^k n_i\Phi_i. $$

To see this, we start with the formula

$$ \Psi\Theta = \Psi( 1_Q\lambda)^G = (\Psi_{QE}(
1_Q\lambda))^G, $$

\noindent which holds because $\Theta$ is induced in \eqref{Theta} from the character $
1_Q\lambda$ of $QE$. Since $\Psi$ is inflated \eqref{Psi} from $\psi \in
\Irr(P)$, its restriction $\Psi_{QE}$ is inflated from the character
$\psi_Q$ of $Q \simeq QE/E$. It follows that $\Psi_{QE}(
1_Q\lambda)$ is precisely the product character $\psi_Q\lambda$. In view
of (i) this product character is

$$\Psi_{QE}(1_Q\lambda) = \psi_Q\lambda = \sum_{i=1}^k
n_i\phi_i\lambda. $$

Hence

$$ \Psi\Theta = (\sum_{i=1}^k n_i\phi_i\lambda)^G = \sum_{i=1}^k
n_i(\phi_i\lambda)^G = \sum_{i=1}^k n_i\Phi_i. $$

Thus (ii) holds.
\end{proof}

\begin{example}\label{example1}
Fix  a prime number $p$. Let $k>0$ and $n_i>0$, $i=1,\ldots, k-1$, be integers.
Choose integers $r, t >0$ such that $p^r > \sum_{i=1}^{k-1} n_i$ and
$p^t\geq k$. Let $Z_{p^t}$ be a cyclic group of order 
$p^t$. Let $N= Z_{p^t}\times Z_{p^t}\times \cdots \times Z_{p^t}$ be the 
direct product of $p^r$ copies of $Z_{p^t}$. Thus $|N|=(p^t)^{p^r}$. Observe that 
$Z_{p^r}$ acts on $N$ by permuting the entries of $N$. Set
$P= Z_{p^{r}}\ltimes N=Z_{p^r} N$ and thus $P=Z_{p^t}\wr Z_{p^r}$ is the wreath product of $Z_{p^t}$
by $Z_{p^r}$.
 
Fix a generator $c$  in $Z_{p^t}$ and  $\alpha\in \Lin(Z_{p^t})$ such that
$\alpha(c)$ is a primitive $p^t$-th root of unity. Set
$\lambda=(\alpha, 1_{Z_{p^t}}, \ldots, 1_{Z_{p^t}})\in \Lin(N)$. 
 We can check that the stabilizer 
$P_{\lambda}$ of $\lambda$ in $P$ is $N$, and so $\psi=\lambda^P\in \Irr(P)$
is a character of degree $p^r$. Observe also that $\psi_N=\sum_{i=1}^{p^r} \lambda_i$,
where $\lambda_i=( 1_{Z_{p^t}}, \ldots, 1_{Z_{p^t}}, \alpha, 1_{Z_{p^t}},\ldots, 1_{Z_{p^t}})\in \Lin(N)$ is the character with
$\alpha$ in the $i$-th position, for $i=1,\ldots, p^r$.

Fix 
\begin{equation*}
q=(\overbrace{c,\ldots, c}^{n_1}, \overbrace{ c^2,\ldots, c^2}^{n_2},\overbrace{ c^3,\ldots, c^3}^{n_3}, \ldots,\overbrace{c^{k-1},\ldots, c^{k-1}}^{n_{k-1}}, 1, 1,\ldots, 1)\in N,
\end{equation*}
\noindent where the first $n_1$-entries are $c$, the $(n_1+1)$-th entry to the $(n_1+n_2)$-th is 
$c^2$ and so for. 

Let $Q$ be the subgroup of $N$ generated by $q$. Observe that $Q$ is a cyclic subgroup 
of $N$ of order $p^t$. Let $\delta\in \Lin(Q)$ such that $\delta(q)=\alpha(c)$.
Since $p^t\geq k$ and $\alpha(c)$ is a primitive $p^t$-root of unity, then for $1\leq i,j\leq k-1$ we have that
$\delta^i\neq \delta^j$ if $i\neq j$. 
 Observe that $\lambda_l(q)=\alpha(c)=\delta(q)$ for $1\leq l\leq n_1$ and so 
$n_1=(\psi_Q, \delta)$. For $l>n_1$, we can check that
  $\lambda_l(q)=\alpha(c^j)=\alpha^j(c)=\delta^j(q)$ if and only if 
$\sum_{i=1}^{j-1} n_i< l\leq \sum_{i=1}^{j} n_i$ and so 
$n_j=(\psi_Q, \delta^j)$. Observe that if
$l> \sum_{i=1}^{k-1} n_i$, then $\lambda_l(q)=\alpha(1)=1$ and so 
$n_k=p^r-\sum_{i=1}^{k-1}n_i= (\psi_Q, 1_Q)$.
Thus $\psi_Q= \sum_{i=1}^{k-1}n_i \delta^i+ n_k 1_Q$. 
\end{example}

Observe that Theorem A follows then by the previous example and Theorem B.
\end{section}

 \end{document}